\documentclass[11pt,reqno]{amsproc} 
 
\usepackage{amstext,amsmath,amssymb,amsthm}

\usepackage{bm}

\usepackage{tikz}
\usepackage{tikz-3dplot}

\usepackage{latexsym}
\usepackage{exscale}
\usepackage{color}

\newcommand{\Primes}{{ \mathbb  P  }}

\newcommand{\C}{  \mathbb  C }
\newcommand{\Z}{  \mathbb Z }
\newcommand{\N}{  \mathbb N }
\newcommand{\Q}{  \mathbb Q }

\newtheorem{Thm}{Theorem}[section]
\newtheorem{Lemma}[Thm]{Lemma}
\newtheorem{Cor}[Thm]{Corollary}

\newcommand{\dsize}{\displaystyle}
\newcommand{\cal}{\mathcal}
\numberwithin{equation}{section}
\newcommand{\E} {{\mathcal E}}

\newcommand{\Mod}[1]{\, (\mathrm{mod}\, #1)}

\theoremstyle{remark}
\newtheorem{Rem}{Remark}[section]

\begin{document}
	\title[Divisibility criteria for cyclotomic polynomials ] {\large  Divisibility criteria and coefficient formulas for cyclotomic polynomials}

	\author{L. De Carli}
	\address{Laura De Carli, Department of Mathematics, Florida International University,   Miami, FL 33199, USA.}
	\email{decarlil@fiu.edu}

	\author{M.   Laporta}
	\address{Maurizio Laporta, Dipartimento di Matematica e Applicazioni, Università degli Studi di Napoli  "Federico II", 80126 Napoli, Italy}
	\email{mlaporta@unina.it }
	
	
		
		\maketitle

		\begin{abstract}
We establish  necessary and sufficient conditions  for a polynomial to be divisible by a cyclotomic polynomials and  derive new formulas involving Ramanujan sums as an application of our results. Additionally, we   provide  new insights into the coefficients of cyclotomic polynomials and we propose a recursive relation between the coefficients of two cyclotomic polynomials whose  indexes differ by a prime factor.

	\end{abstract}

	\section{Introduction}
	The {\it cyclotomic polynomial} of {\it index} $n\in\N$  is defined 
	as
	$$
	\Phi_n(z):=\prod_{1\leq j\leq n\atop{(j,n)=1}} (z-\zeta_n^j),
	$$
	where $(j,n)$ denotes the greatest common divisor of   $j,n\in\N$  and $\zeta_n^j:=e^{2\pi i j/n}$ is a $n$th primitive root of unity  for  $(j,n)=1$.  
	The degree of $\Phi_n$ is given by the {\it Euler totient function}	$\varphi(n):=\#\{j\in\N\cap[1,n]:\ (j,n)=1\}$, and we write
	\begin{equation}\label{explicitcyclotomic}
	\Phi_n(z)=\sum _{k=0}^{\varphi(n)}a_n(k)z^{\varphi(n)-k},
	\end{equation}
	with $a_n(0)=1$.  It is well-known that $\Phi_n(z)$ is irreducible  over $\Q$ and $a_n(k)\in\Z$ for any $k$.

	 The {\it order} of $\Phi_n$ is  the number of distinct prime factors of $n$, which is denoted by $\omega(n)$. If $n$ is square-free,   then $\Phi_n$ is referred to as {\it binary, ternary}, etc., when $\omega(n) = 2, 3$, and so on, respectively.

	The study of cyclotomic polynomials has a very long history, which goes back at least to Gauss. We  refer the reader to the surveys by  C. Sanna \cite{S} and R. Thangadurai \cite{T}  for an overview of  results on these polynomials.  We use here two basic properties which are stated
	in Lemma \ref{cycloformulas} below.
	
	 \subsection{Divisibility results}
	   Understanding when a given polynomial is divisible by a cyclotomic polynomial $\Phi_n(x)$ is equivalent to understanding when it vanishes at all primitive $n$-th roots of unity, an important condition that arises in a wide range of pure and applied mathematical problems.  
	 
	  For instance,  conditions of  divisibility by   cyclotomic polynomials   are  closely related to a conjecture on tiling of the integers proposed by Coven and Meyerowitz  \cite{CM} and to the  {\it Fuglede’s spectral set  conjecture}  in dimension $d=1$ \cite{F}.  For   details on these conjectures  and the  connections between  them,  see  \cite{CM, L} and the references cited there.

	 In  this paper we provide new necessary and sufficient conditions for a   polynomial with complex coefficients,
	 \begin{equation}\label{polynomial}
	 	{\cal P}(z):=  a_0z^m+a_1z^{m-1}+\ldots+a_{m-1}z+a_m,
	 \end{equation} 
	 to be divisible by a given cyclotomic polynomial.
	 
	  Our main result is the following
	  
	  \begin{Thm}\label{T-conjecture} Let $ N\in\N$ be such that
	  	 $  \varphi(N)\le m={\deg}{\cal P}(z)$.	
	  	 
	  	The polynomial  ${\cal P}(z) $ is divisible by the cyclotomic polynomial $
	  	\Phi_N(z )$ 
	  	if and only if 
	  	\begin{equation}\label{e-conj-p>2}
	  		\sum_{d|N}\mu(d)\sum_{j\equiv_{N }\ h-\sum_{p|d}\frac{N }{p}} a_j=0
	  	\end{equation}
	  	for every $h\in  \{0, 1,\ldots,  N -1\}$.

	  	When $N$ is even,    \eqref{e-conj-p>2} needs to be satisfied only for every $h\in \{0,\ldots, N/2 -1\}$.
	  	
	  \end{Thm}
	  
	 Let us give some clarifications about the notation in \eqref{e-conj-p>2} and the rest of the paper.

	 As usual in Number Theory, the letter $p$ (with or without subscript)  is reserved for the prime numbers. We denote the set of primes numbers  by $\Primes$.

	 	In \eqref{polynomial}, \eqref{e-conj-p>2} and in what follows, it is understood that $a_j=0$ when $j>m$ or $j<0$. Further, we adopt the convention that
	 $\sum_{p|d}N/p=0$ when $d=1$.

	 Recall that $\mu$ is  the  {\it M\"obius function} defined as
	$\mu(n)=(-1)^{\omega(n)}$ if $n$ is square-free, $\mu(n)=0$ otherwise, where $\omega(1)=0$ and
	$\omega(n)$ is
	the number of the distinct prime factors of $n\ge 2$. In particular, we see that the equation \eqref{e-conj-p>2} can be equivalently written as
	$$
	\sum_{d|\gamma(N)}(-1)^{\omega(d)}\sum_{j\equiv_{N }\ h-\sum_{p|d}\frac{N }{p}} a_j=0,
	$$
	where $$
	\gamma (N):=
		\begin{cases} 1& \mbox{if $N=1$},\cr
		\prod_{p|N}p  & \mbox{if $N>1$,}
		\cr\end{cases}
		$$ 
		is the so-called {\it (square-free) kernel} of $N$.

For brevity, we often write  $m\equiv_k n$ or $m\equiv n\ (k)$  to mean that  $m\equiv n \Mod k$, i.e.,
	 $k$ divides $m-n$.  If, in addition, $0\le n<k$, we also write $n=\{m\}_k$, i.e.,   $$m=k[m/k] +\{m\}_k,$$ where $[x]$ denotes the integer part of the real number $x$.
	\medskip
 
As a consequence of Theorem~\ref{T-conjecture}, we derive new identities involving the so-called {\it Ramanujan sums} \cite{R}
\begin{equation}\label{def-Ram-sums}
	c_n(r):=\sum_{\substack{1\le j \le n\atop{(j,n)=1}}} \zeta_n^{jr}, \quad (n\in\N, \, r\in\Z).
\end{equation}
Some of their basic properties are summarized in Lemma~\ref{rama} below.  

Indeed, we exploit a notable result of Tóth  \cite{To}, who has shown  that the polynomial
\begin{equation}\label{Toth1}	
	T_n(z):=\sum_{r=0}^{n-1}c_n(r)z^r - n\end{equation} 
is divisible by  the cyclotomic  polynomial 
$\Phi_n(z)$. The degree of $T_n(z)$ is
$\tau =\tau(n):= n-\frac{n}{\gamma(n)}$  \cite[Theorem 3]{To}, where $\gamma (n)$ is the kernel of $n$. Thus, by applying 
Theorem~\ref{T-conjecture} to $T_n(z)$ we get the following

  \begin{Cor}\label{T-Ram-conj}
  	Given $n\in\N$,  for every divisor $d_1$ of $\gamma(n)$ we have that	
  	\begin{equation}\label{directidentity}		
  		\mu(d_1)\sum_{d|n\atop d\not=d_1}\mu(d)c_n\Big(\sum_{p|d_1}\frac{n}{p}-\sum_{p|d}\frac{n}{p}\Big)=n-\varphi(n).
  	\end{equation}
  	Further, if
  	$$
  	h\in H=H(n):=
  	\begin{cases} \{0, \ldots, n/2-1\} & \mbox{if $n$ is even},\cr
  		\{0, \ldots, n-1\}   & \mbox{otherwise,}
  		\cr\end{cases}
  	$$
  	 is  such that $\{h- \sum_{p|d}n/p\}_n < n-n/\gamma(n)$ 
  	 for every divisor $d$ of $\gamma(n)$,   then
  	
  	\begin{equation}\label{e-sum-cn}
  		\sum_{d|n}\mu(d)  c_n\Big(h+\frac n{\gamma(n)}  -\sum_{p|d}\frac{n}{p}\ \Big)=0.
  	\end{equation} 
  If, in addition, $\mu(n)=0$, then for every divisor $d_1>1$ of $\gamma(n)$ we have that	
  	\begin{equation}\label{e-sum-cnnonsquarefree}
  	\sum_{d|n\atop d\not=d_1}\mu(d)  c_n\Big(h+\frac n{\gamma(n)}  -\sum_{p|d}\frac{n}{p}\ \Big)=0
  \end{equation}
  for all $  h\in\big(\sum_{p|d_1}n/p-n/\gamma(n),\ \sum_{p|d_1}n/p\big)\cap H$.
  \end{Cor}

  	In Section \ref{appendix} we give a longer proof of
  	\eqref{directidentity}  by using only the basic properties of the Ramanujan sums and the functions $\mu$, $\varphi$.
  	However, we think that \eqref{e-sum-cn} and \eqref{e-sum-cnnonsquarefree} cannot be established in the same elementary fashion.

\subsection{Results on the coefficients of cyclotomic polynomials}
  The coefficients of cyclotomic polynomials   are the subject of intensive study and many formulas are
known for them.   We refer the reader to the survey by Herrera-Poyatos and Moree \cite{HPM}, which includes   most of the formulas where
Bernoulli numbers, Stirling numbers, and Ramanujan's sums are involved. 
Recursive relations between the coefficients of two cyclotomic polynomials  are key 
to the so-called big prime algorithm of Arnold and Monagan \cite{AM}.

In this paper  we propose new formulas for the  coefficients of cyclotomic polynomials. We also  establish an alternate version of the Arnold and Monagan formula (see \eqref{e-formula-AM} below) by using a generalization of Vieta's formulas introduced in~\cite{DEL}  in terms of the  {\it complete homogeneous symmetric  polynomials} ${\cal H}_r(z_1,\ldots, z_n)$, with  $n\in\N$ and $r\in\N_0:=\N\cup\{0\}$.   
We have used  \cite[\S 1.2]{Mcdonald} for the definition and the main properties of  these polynomials. 

Recall that
${\cal H}_r(z_1,\ldots, z_n)$, for $r\in\N$,
is the sum of all monomials of   degree $r$ in the variables $z_1,\ldots, z_n$. We let 
$$
{\cal H}_r (z_1,\ldots,z_n) :=  \dsize\sum_{1\le j_1\leq j_2\leq \ldots\leq j_r\le n}z_{j_1} z_{j_2} \cdots z_{j_r} 
= \dsize\sum_{r_1+r_2\ldots+ r_n=r\atop r_i\in\N_0}z_1^{r_1}z_2^{r_2}\cdots z_n^{r_n}.
$$
When $r=0$, it is  ${\cal H}_0(z_1,\ldots,z_n):=1$ for all $z_1,\ldots, z_n$. 
Note that ${\cal H}_{r}(z)=z^r$ for each $r\in\N_0$.

Let us also recall that the {\it elementary symmetric polynomials} of degree $m\in\N_0$ are defined as
$$
{\cal E}_m(z_1,\ldots,z_n):= \sum_{1\le j_1<j_2<\ldots< j_m\le n}z_{j_1}\cdots z_{j_m}
$$
for $1\le m \leq n$,
while $ {\cal E}_0(z_1,\ldots, z_n):=1$ for all $z_1,\ldots, z_n$. Further, we let ${\cal E}_k(z_1,\ldots, z_n):=0$ whenever  $k>n$.

Both  ${\cal H}_m (z_1,\ldots,z_n)$ and ${\cal E}_m (z_1,\ldots,z_n)$  are  homogeneous of degree $m$ and invariant under permutation of the variables $z_j$.    They are related by   the identity
$$
\sum_{j=0}^m(-1)^j{\cal E}_j(z_1,\ldots, z_n){\cal H}_{m-j}(z_1,\ldots, z_n)=0.
$$
We use  these properties throughout the paper without mentioning them explicitly.  

It is well-known that the coefficients of any 
polynomial can be expressed in terms of its roots by means of the elementary symmetric polynomials.
Specifically, if \ $z_1,\ldots, z_{m}$ are the roots (not necessarily distinct) of  the polynomial ${\cal P}(z)$  given in \eqref{polynomial},
then  

\begin{equation}
	(-1)^k\frac{a_k}{a_0}= {\cal E}_k(z_1,\ldots, z_{m}),\quad \forall k\in\{0,1,\ldots,m\}.\label{e-Vieta}
\end{equation}
(see e.g. \cite[Ex. 4.6.6]{B}).
This formula is usually named after Fran\c{c}ois Vi\`ete (1540-1603)  more commonly referred to by the Latinized form of his name Franciscus Vieta.

 In \cite{DEL} we have  generalized Vieta's formula  by establishing
necessary conditions that the coefficients of ${\cal P}(z)$  need to satisfy when only some of the roots of ${\cal P}(z)$ are given. 
See Lemma \ref{L-division-pol-3} below for such a generalization.

\medskip

Before presenting our results on the coefficients of  cyclotomic polynomials, we introduce some   notation.

For $\zeta_n  :=e^{2\pi i /n}$, let us  define the vectors 

\begin{itemize}
	\item
	$\vec{\zeta}_n:=(1,\zeta_n,\ldots ,\zeta_n^j,\ldots, \zeta_n^{n-1})$
	
	\item
	$\vec\zeta_n^*:=(\zeta_{1,n}^*,\ldots ,\zeta_{j,n}^*,\ldots,\zeta_{n,n}^*)$, where    
	$ \dsize
	\zeta_{j,n}^*:=\begin{cases} \zeta_n^j  & \mbox{\ \hbox{if}\ $(j,n)=1$},\cr
		0& \mbox{otherwise.}
		\cr\end{cases}
	$ 
	\item $\vec\zeta_n^{**}:=(\zeta_{1,n}^{**},\ldots ,\zeta_{j,n}^{**},\ldots,\zeta_{n,n}^{**})$, where
	$ \dsize
	\zeta_{j,n}^{**}:=\begin{cases} \zeta_n^j  & \mbox{if $(j,n)>1$},\cr
		0& \mbox{otherwise.}
		\cr\end{cases}
	$ 	
\end{itemize}
Clearly, $\vec\zeta_n=\vec\zeta_n^{*}+\vec\zeta_n^{**}$.

\medskip
For any given function $f:\C^n\to\C$, 
we write $f(\vec v)=f(v_1,\ldots,v_n)$ instead of $f(\vec v)=f\big((v_1,\ldots,v_n)\big)$.

  {Since ${\cal H}_k(x_1,x_2,\ldots,x_h,0,\ldots,0)={\cal H}_k(x_1,x_2,\ldots,x_h)$, it turns out that
	${\cal H}_k(\vec\zeta_n^*)$ is evaluated only at the $n$th primitive roots of unity. 
	Analogous considerations hold  for ${\cal E}_k(\vec\zeta_n^*)$,  ${\cal H}_k(\vec\zeta_n^{**})$, and ${\cal E}_k(\vec\zeta_n^{**})$.}
\medskip

From \eqref{e-Vieta} and our generalization of  Vieta's formula applied to the cyclotomic polynomial  \eqref{explicitcyclotomic}, it follows  immediately 
that
\begin{equation}\label{coeffNewton}
	a_n(k)=(-1)^{k}{\cal E}_{k}(\vec\zeta_n^*)={\cal H}_k(\vec\zeta_n^{**})
\end{equation}  
for every $k\in\{0,\ldots,\varphi(n)\}$.   {See Lemma \ref{L-rootsformulas}}.
Since  the coefficients of the cyclotomic polynomials are integers,  the formula \eqref{coeffNewton}   implies that so are 
${\cal E}_{k}(\vec\zeta_n^*), {\cal H}_k(\vec\zeta_n^{**}) $ for every $n\in\N$ and every $k\in\{0,\ldots,\varphi(n)\}$.
\medskip

We can now  state our recursive relation between the coefficients of two cyclotomic polynomials whose  indexes differ by a prime factor.

\begin{Thm}\label{T-newcoeffformula}  Assume that $p\in\Primes$
	does not divide $m\in\N$.  The coefficients of the cyclotomic polynomial 
	$$
	\Phi_{mp}(z):=\sum _{k=0}^{(p-1)\varphi(m)}a_{mp}(k)z^{(p-1)\varphi(m)-k}=\prod_{1\leq j\leq mp\atop{(j,mp)=1}} (z-\zeta_{mp}^j)
	$$
	are given by
	\begin{equation}\label{newAM}
		a_{mp}(k)=
		\sum_{s=0}^{[k/p]} a_m(s){\cal H}_{k-sp}(\vec\zeta_m^*),
	\end{equation}
	with
	$k\in\{0,1,\ldots, (p-1)\varphi(m)\}$.
\end{Thm}

\begin{Rem}  Under the same hypothesis of Theorem \ref{T-newcoeffformula}, the aforementioned recursive formula of Arnold and Monagan  is given in the form 
	\begin{equation}\label{e-formula-AM}
		a_{mp}(k)=	a_{mp}(k-m)-
		\sum_{j,h\ge 0\atop jp+h=k}a_m(j)b_{m}(h),
	\end{equation}
	where  the $b_{m}(0),b_{m}(1),\ldots,b_{m}(m-\varphi(m))$ are the coefficients
	of the so-called   $m$th inverse cyclotomic polynomial
	$
	\Psi_m(z):=(z^m-1)/\Phi_m(z)$. See \cite[\S 4]{AM}.
	Both our  proof of   \eqref{newAM}  and   the proof of    \eqref{e-formula-AM} by Arnold and Monagan start  with the identity $\Phi_{m}(z^p)=\Phi_{mp}(z)\Phi_m(z)$ (see Lemma \ref{cycloformulas} below),  but    then  Arnold and Monagan proceed   with the identity
	$$
	\Phi_{mp}(z)=\frac{\Phi_{m}(z^p)}{\Phi_m(z)}=
	\Phi_{m}(z^p)\frac{\Psi_m(z)}{z^m-1}=
	-\Phi_{m}(z^p)\Psi_m(z)\sum_{s\ge 0}z^{sm},
	$$
	and instead we use 
	our generalization of Vieta's formulas (Lemma \ref{L-division-pol-3}).    
	
	Further, from \eqref{coeffNewton} and  
	\eqref{newAM} we 
	obtain the following recursive formula involving the roots of unity when $p$ does not divide $m$:
	$$
	{\cal H}_k(\vec\zeta_{mp}^{**})=
	\sum_{s=0}^{[k/p]}{\cal H}_{s}(\vec\zeta_m^{**}){\cal H}_{k-sp}(\vec\zeta_m^*)
	$$
	for every 
	$k\in\{0,1,\ldots,  (p-1)\varphi(m)\}$.
\end{Rem}
\smallskip

\begin{Rem}
	If  $p$ divides $m$, then $\Phi_{mp}(z)=\Phi_{m}(z^p)$ (see
	Lemma \ref{cycloformulas}), and  ${\rm deg}\Phi_{mp}=\varphi(mp)=\varphi(m)\varphi(p)\frac{(m,p)}{\varphi(m,p)}=
	p\varphi(m)$ (see Lemma \ref{eulerandmebius}). Consequently,
	$$
	\Phi_{mp}(z)=\sum _{k=0}^{p\varphi(m)}a_{mp}(k)z^{p\varphi(m)-k}=	\Phi_{m}(z^p)=\sum _{s=0}^{\varphi(m)}a_m(s)z^{p\varphi(m)-ps}
	$$
	from which it follows immediately that
	$$
	a_{mp}(k)
	=\begin{cases} a_m(s)& \mbox{if $k=ps$,}\cr
		0& \mbox{otherwise.}
	\end{cases}
	$$ 
	In view of  \eqref{coeffNewton}, when $p$ divides $m$  the following formula holds:
	$$
	{\cal H}_k(\vec\zeta_{mp}^{**})
	=\begin{cases} {\cal H}_{s}(\vec\zeta_{m}^{**}) & \mbox{if $k=ps$},\cr
		0& \mbox{otherwise,}
	\end{cases}
	$$
	for every 
	$k\in\{0,1,\ldots,  p\varphi(m)\}$.
\end{Rem}
\smallskip

The special case of Theorem \ref{T-newcoeffformula}, with $m\in\Primes$, concerns the binary cyclotomic polynomials.
It is well known that the coefficients of such polynomials lie in the set $\{-1,\, 0,\, 1\}$. This was first proved by Migotti \cite{Mi} in 1883, and has since been reproved and extended by various authors; see  for example  \cite{Beiter, LL}. Here we exhibit another simple and direct proof of this result.

\begin{Thm}\label{T-coeff-pq}
	Let $p,q\in\Primes$, with $p\not=q$. 
	
	For every $k\in\{0,1,\ldots,\varphi(pq)\}$ one has
	$$
		a_{pq}(k)
		={\cal H}_k(1, \zeta_p,\ldots, \zeta_p^{p-1}, \zeta_q,\ldots, \zeta_q^{q-1})\in\{-1,\, 0,\, 1\}.
	$$
\end{Thm}

Our  paper is organized as follows. In Section \ref{lemmata} we  provide  some lemmas that are essential for our proofs. In Section \ref{proofs} our main results are proved.
In Section \ref{appendix} we give both an alternate proof of 
\eqref{directidentity} and the instance of Corollary \ref{T-Ram-conj} for $\omega(n)=2$.

\medskip
In closing, we would like to point out that  we  use   $(m,n)$ to denote  the greatest common divisor of the integers $m$ and $n$, but   also to   denote an  open interval with endpoints $x,y$, or a vector with components $x,y$. The meaning will always be evident from the context.

\section{Lemmata}\label{lemmata}

First, we recall some basic properties of  the  M\"obius function $\mu$ and the Euler totient function $\varphi$. For the proof see \cite{A}.  Most of the properties recalled in the next two lemmas are often used here
without mentioning them explicitly. 

\begin{Lemma}\label{eulerandmebius} 
	  The following identities hold for any $n\in\N $.
		
	\begin{equation}
	\sum_{d|n}\mu(d)=
	\begin{cases}
		1 &\text{if $n=1$,}\\
		0 &\text{otherwise.}
	\end{cases}\label{mertens}
		\end{equation}
	
	\begin{equation}
		\frac{\varphi(n)}{n}=\sum_{d|n}\frac{\mu(d)}{d}=\prod_{p|n\atop p\in\Primes}\Big(1-\frac{1}{p}\Big).	\label{totient}
	\end{equation}

		\begin{equation}
	\frac{n}{\varphi(n)}=\sum_{d|n}\frac{\mu(d)^2}{\varphi(d)}.\label{squarefree}
		\end{equation}

\end{Lemma}

In particular, \eqref{totient} implies that 
$$
\varphi(mn)=\varphi(m)\varphi(n)\frac{(m,n)}{\varphi(m,n)}
$$
for all $m,n\in\N$. Thus, 
$\varphi$ is a multiplicative function, i.e., $\varphi(mn)=\varphi(m)\varphi(n)$ when $(m,n)=1$. It is readily seen by its definition that $\mu$ is a multiplicative function, as well.
\smallskip

Now, let us recall some basic properties of the Ramanujan sums \eqref{def-Ram-sums}.  For a comprehensive treatise on them, we refer the reader
to \cite{MC}.

\begin{Lemma}\label{rama} Let $n,m\in\N$, $h,r\in\Z$. We have that
	\begin{itemize}
		\item[(1)] $c_n(0)=\varphi(n)$
		\item[(2)] $c_n(r)=c_n(-r)$
		\item[(3)] $c_n(r+hn)=c_n(r)$

		\item[(4)] (H\"older's identity)\quad $\dsize c_n(r)=\varphi(n)\frac{\mu(n/d)}{\varphi(n/d)}$, where $d:=(n,r)$
		\item[(5)] If $(n,r)=1$, then $c_n(r)=\mu(n)$
		\item[(6)] If $(n, m)=1$, then $c_{nm}(r)=c_n(r)c_m(r)$ 
		\item[(7)] $c_n\big(n-n/\gamma(n)\big)=c_n\big(n/\gamma(n)\big)=(-1)^{\omega(n)}\frac{n}{\gamma(n)}$, where 
		$\gamma(n)$ 
		is the kernel of $n$.
		\item[(8)] If  $n$ is even, then 
		$c_n(h+mn/2)=(-1)^mc_n(h)$.
	\end{itemize}
\end{Lemma}

\begin{proof} For the properties (1)-(6) we refer the reader to \cite{A} and \cite{MC}. The first equation of (7)
	follows immediately from (2) and (3). The second equation is established by applying H\"older's identity (4) with 
	$r=n/\gamma(n)$, so that $d:=(n,r)=n/\gamma(n)$ and
	$$
	c_n\big(n/\gamma(n)\big)=
	\varphi(n)\frac{\mu\big(\gamma(n)\big)}{\varphi\big(\gamma(n)\big)}.
	$$
	The conclusion follows 
	after noticing that \eqref{totient} yields
	$\dsize
	\frac{\varphi(n)}{\varphi\big(\gamma(n)\big)}=\frac{n}{\gamma(n)}$.
	
	Finally,  in order to prove (8) it suffices to observe that 
	$$
	c_n\Big(h+\frac{mn}{2}\Big):=\sum_{j=1\atop (j,n)=1}^n\zeta_n^{j(h+\frac{mn}{2})}
	=\sum_{j=1\atop (j,n)=1}^n\zeta_n^{jh}e^{\pi ijm},
	$$
	where $e^{\pi i jm}=(-1)^m$ because $j$ is odd. 
\end{proof}

In the next lemma we state  two  well-known properties of the  cyclotomic polynomials. For the proof see \cite{S, T}.

\begin{Lemma}\label{cycloformulas} Let  $m, n, k\in\N$ and $p\in\Primes$.  We have
	\begin{equation}\label{e-with moebius}
		\Phi_n(z)= \prod_{d\vert n} (z^{n/d}-1)^{\mu(d)}.
	\end{equation}
	
		Further,
	\begin{equation}\label {e-primepower}
		\Phi_{mp^k }(z)=\begin{cases} \Phi_m(z^{p^{k }})  & \mbox{if $p|m$},\cr\cr
			\dsize	\frac{\Phi_m(z^{p^k})}{\Phi_m(z^{p^{k-1}})}  & \mbox{otherwise.}
			\cr\end{cases}
	\end{equation}
	
\end{Lemma}

\begin{Rem} As a consequence of \eqref{e-with moebius}, for $p\in\Primes$ one has 
	\begin{equation}\label {e-p-prime}
		\Phi_p(z)=\frac{z^p-1}{\Phi_1(z)}=\frac{z^p-1}{z-1}=1+z+\ldots+ z^{p-1}.
	\end{equation}

	Moreover, for all $m,n\in\N$, 	we  have that
	\begin{equation}\label{e-p-primes}
		\Phi_{mn }(z)= \Phi_{m\gamma(n)  }(z^{n/\gamma(n)} ),
	\end{equation}
	where $\gamma(n)$ is the kernel of $n$.
	Indeed, recalling that $\mu(d)=0$ if $d$ is not square-free and using  \eqref{e-with moebius}, we see that
		$$
	\Phi_{mn }(z)=\prod_{d\vert m\gamma(n)} (z^{mn/d}-1)^{\mu(d)}.
	$$
	Thus, \eqref{e-p-primes} follows after writing 
	$ 
	z^{mn/d}=(z^{n/\gamma(n)})^{m\gamma(n)/d}
	$ 
	and applying \eqref{e-with moebius} again.

	In particular, the  identity \eqref{e-p-primes}  shows that in order to study the coefficients of cyclotomic polynomials it suffices to consider  those with square-free  index. 
\end{Rem} 
\smallskip

The  following properties  of the complete homogeneous symmetric polynomials
are proved   in \cite{K}.

\begin{Lemma}\label{L-id-homog}\indent
	
	\begin{itemize}
		\item Given $n\in\N$ we have that
		\begin{equation}\label{e-id-2terms} {\cal H}_{n}(x_1,\ldots,x_m)= x_1{\cal H}_{n-1}(x_1,\ldots,x_m) + {\cal H}_{n }(x_2,\ldots,x_m).
		\end{equation}
		In particular, for $x_1=1$ and any $s\in\N\cap[1,n]$ this yields
		$$
		{\cal H}_{n }(1,x_2,\ldots,x_m)-{\cal H}_{n-s}(1,x_2,\ldots,x_m) =\sum_{k=0}^{s-1} {\cal H}_{n-k}(x_2,\ldots,x_m).
		$$
		\item 
		For $n\in\N_0$ we have 
		\begin{equation}\label{e-sigma-id 2}  \dsize
			\hskip - 1 mm	{\cal H}_{n}(x_1,\ldots,x_{m_1},\, y_1,\ldots,y_{m_2})= \sum_{s=0}^n{\cal H}_s (x_1,\ldots,x_{m_1}) {\cal H}_{n-s}(y_1,\ldots,y_{m_2}).
		\end{equation}
	\end{itemize}
\end{Lemma}
The  next three lemmas are   proved  in \cite{DEL}.

\begin{Lemma}\label{L-1.5}
	Let $n\in\N $, $n\ge 2$. Given $k\in\N$ and $m\in\N_0$,   
	$$
	{\cal H}_m(1, \zeta_n,\ldots, \zeta_n^{k})
	=\begin{cases}  1 &  \mbox{if either $\{k\}_n=0$ or $ \{m\}_n =0$,}   
	\cr\dsize\prod_{s=1}^{\{k\}_n  }\frac{1-\zeta_n^{\{m\}_n +s }}{1-\zeta_n^s } &\mbox{if   $1\leq  \{k\}_n +\{m\}_n< n $,}\cr  0&   \mbox{if $ \{k\}_n +\{m\}_n \ge n$. }\cr\end{cases}
$$
\end{Lemma}
Recall that  $\{k\}_n :=k-n[k/n]$.
 
The particular instance $k=n-1$ gives 
\begin{equation}\label{e-sigma-id 3} {\cal H}_m(\vec\zeta_n) =\begin{cases}  1 &  \mbox{if   $n|m$,}
	\cr  0&   \mbox{otherwise. }\cr\end{cases}\end{equation}

The next lemma is \cite[Th.\ 1.1 ]{DEL}. It generalizes Vieta's formula  by providing  necessary conditions
that the coefficients of  the polynomial ${\cal P}(z) $ in \eqref{polynomial} 
need to satisfy when
 some of its roots are given.

\begin{Lemma}\label{L-division-pol-3}
If  $z_0,\ldots, z_n\in \C$ are roots of 
$ 
{\cal P}(z) $, 
then for any $s\in\{0,1,\ldots, n\}$ one has
 
$$ 
 \sum_{j=0}^{m-s } a_{j}  {\cal H}_{m-s-j}(z_0,\ldots, z_{s})=0.
$$ 
Further,
$$ 
{\cal P}(z)=  \prod_{k=0}^n
(z-z_k)\sum_{j=0}^{m-n-1}c_jz^{m-n-j},
$$
where 
\begin{equation}\label{e-ck1}    c_{k}=    \sum_{j=0}^k a_{ j}{\cal H}_{k-j}(z _0,\ldots, z_n)
\end{equation}
for any $k\in\{0,\ldots,m-n-1\}$.

\end{Lemma}

The following  lemma is \cite[Corollary 1.2]{DEL}.   	
\begin{Lemma}\label{L-rootsformulas} Let $n\in\N$.
	
	\begin{itemize}
		\item  If $0\le k\leq \varphi(n)$, then
		$${\cal E}_{k}(\vec\zeta_n^*)=(-1)^{k}	{\cal H}_k(\vec\zeta_n^{**}).$$
		\item  If $0\le k\leq n-\varphi(n)$, then  
		$${\cal E}_{k}(\vec\zeta_n^{**})=(-1)^{k}
		{\cal H}_k(\vec\zeta_n^*).
		$$ 
	\end{itemize}
	
\end{Lemma}

\begin{Rem}
	Both formulas above are  valid only when $k$ belongs to the assigned range.  For example, if  $n$ is prime, then $\E_n(\vec\zeta_n ^*)=0$ by definition,  while ${\cal H}_n(\vec\zeta_n^{**})= {\cal H}_n(1)=1.$

\end{Rem}
\smallskip

We use Lemma \ref{L-division-pol-3} to prove a condition on the divisibility of the polynomial \eqref{polynomial} by $z^n-1$. 

\begin{Lemma} \label{L-cyclic-2} Let $n\in\N$ such that $1\leq n\le m$.

The polynomial 
$z^n-1$ divides ${\cal P}(z)$   if and only if
\begin{equation}\label{id-11} 
	\sum_{ j\equiv r \, (n)} a_j=0
\end{equation}
for every $r\in\{0, 1,\ldots,n-1\}$.
\end{Lemma}

\begin{proof} First, note that $z^n-1$ divides ${\cal P}(z)$   if and only if $1,\zeta_n,\ldots,\zeta_n^{n-1}$ are roots of ${\cal P}(z)$. If this is the case, then for every $h\in\{0, 1,\ldots,n-1\}$	we can  write
\begin{equation*}
	0={\cal P}(\zeta_n^h)=\zeta_n^{hm}\sum_{j=0}^m a_{ j}\zeta_n^{-hj}
	=	\zeta_n^{hm}\sum_{r=0}^{n-1}\zeta_n^{-hr}\alpha_r,
\end{equation*}
where 
$$
\alpha_r:=
\sum_{j=0\atop j\equiv r \, (n)}^m a_j,\quad r=0, 1,\ldots,n-1.
$$
Thus, $z^n-1$ divides ${\cal P}(z)$   if and only if  for each $h=0, 1,\ldots,n-1$  the vectors 
$(1,  \zeta_n^h,\ldots, \zeta_n^{h(n-1)})$ and	$(\alpha_0,\, \alpha_1,\ldots,\alpha_{n-1})$ are  orthogonal.
Since $ \zeta_n^j\ne  \zeta_n^k$ for all distinct integers $j, k\in[0,n-1]$,
the Vandermonde matrix with rows 
$(1,  \zeta_n^h,\ldots, \zeta_n^{h(n-1)})$, $h=0, 1,\ldots,n-1$, is non-singular and so the equation
$$
\left(\begin{matrix} 1&1&... & 1
	\cr 	  1 & \zeta_n& ... & \zeta_n^{n-1}\cr 
	1 &  \zeta_n^2& ... & \zeta_n^{2(n-1)}
	\cr \vdots&\vdots&\vdots & \vdots\cr 1 &  \zeta_n^{n-1}& ... & \zeta_n^{(n-1)^2} 
\end{matrix} \right)	\left(\begin{matrix}\alpha_0&\cr \alpha_1&\cr\alpha_2&\cr\vdots&\cr\alpha_{n-1}\end{matrix} \right)=0
$$
is satisfied only by $	(\alpha_0,\, \alpha_1,\ldots,\alpha_{n-1})=(0,\ldots,0)$, which gives \eqref{id-11}.
\end{proof}

\begin{Rem} 
Using  Lemma \ref{L-cyclic-2},  we can easily prove that  the polynomial
$z^n-\eta$, with $\eta\in\C$, divides ${\cal P}(z)$   if and only if 
\begin{equation}\label{id-111} 
	\sum_{ j\equiv r\, (n)} a_j \eta^{-\frac jn}
\end{equation}
for every $r\in\{0, 1,\ldots,n-1\}$.

Indeed, for  $\eta = |\eta| e^{2\pi i \theta}$, with $0 \leq \theta < 2\pi$, the roots of $z^n - \eta$ are 
$$
\xi_k := \eta^{1/n} \zeta_n^k, \quad k = 0, 1, \ldots, n-1,
$$ 
where $\eta^{1/n} = |\eta|^{1/n} e^{2\pi i \theta/n}$.  
Therefore, the polynomial $z^n - \eta$ divides $\mathcal{P}(z)$ if and only if for every $k\in \{0,1, \ldots, n-1\}$ we have that    
$$ 0= {\cal P}(\xi_k)=
\sum_{j=0}^m a_j\eta^{(m-j)/n} \zeta_n^{(m-j)k}  =\eta^{m/n} 
\, \sum_{j=0}^m b_j \zeta_n^{(m-j)k},
$$ 
where $b_j := a_j\eta^{ -\frac { j}n}$. This implies that that the polynomial  $  \sum_{j=0}^m b_j z^{m-j}$ is divisible by $z^n-1$. 
Lemma \ref{L-cyclic-2} yields
\eqref{id-111}.
\end{Rem}
\smallskip

The following lemma is  key to prove Theorem \ref{T-conjecture}.

\begin{Lemma}\label{T-generaldivisibility}
	Let   $s\in\N$ and $p\in\Primes $ such that $(p-1)s\le m={\deg}{\cal P}(z)$.  
	The polynomial 
	$\Phi_p(z^s)$  divides ${\cal P}(z) $  if and only if 
	\begin{equation}\label{e-cond-oneprime}
		\sum_{ j\equiv_{ps}\  h }a_j=\sum_{ j\equiv_{ps}\  h-s  \ }a_{j}
	\end{equation}
	for every $	h\in\{0, 1,\ldots,  ps-1\}$.
	When  $p=2$, the condition \eqref{e-cond-oneprime} needs to be satisfied only for every $	h\in\{0, 1,\ldots,  s-1\}$.
\end{Lemma}

\begin{proof} 
	For any $p\in\Primes$, the identity \eqref{e-p-prime} yields
	$$ 
	\Phi_p(z^s)= \frac{z^{ps}-1}{z^s-1}.
	$$  
	Thus,  $\Phi_p(z^s)$ divides ${\cal P}(z)$ 	if and only if $z^{ps}-1=\Phi_p(z^s)(z^s-1)$ divides 
	$$
	(z^s-1){\cal P}(z)=\sum_{j=0}^{m+s} b_jz^{m+s-j},	
	$$	
	where
	$b_j:= a_j-a_{j-s}$ for every $j\in\{0,1,\ldots,m+s\}$.
	Recall that we have assumed  $a_r=0$ when $r<0$ or  $r>m$.  
	
	In view of  Lemma \ref{L-cyclic-2} we see that $z^{ps}-1$ divides $(z^s-1){\cal P}(z)$  if and only if 
	$$
	\sum_{j=0\atop j\equiv h \ (ps)}^{m+s} b_j=
	\sum_{j=0\atop j\equiv h \ (ps)}^{m+s}(a_j-a_{j-s})=0
	$$
	for every $h\in\{0, 1, \ldots, ps-1\}$. Thus, 
	$$
	\sum_{j=0\atop j\equiv h \ (ps)}^{m}\!a_j =
	\sum_{j=0\atop j\equiv h \ (ps)}^{m+s}\!a_j=
	\sum_{j=0\atop j\equiv h \ (ps)}^{m+s}\!a_{j-s}\\
	=
	\sum_{j=-s\atop j\equiv h-s \ (ps)}^{m}a_{j}	=\sum_{j=0\atop j\equiv h-s \ (ps)}^{m}a_{j},		
	$$
	for every $h\in\{0, 1, \ldots, ps-1\}$, which is
	 \eqref{e-cond-oneprime}.

	When $p = 2$, the previous condition becomes
	$$
	\sum_{j=0\atop j \equiv h \ (2s)}^m a_j = \sum_{j=0\atop j \equiv h -s \ (2s)}^m a_j, 
	$$
	for every $h\in\{0, 1, \ldots, 2s-1\}$. Let us show that in this case it suffices to take just 
	$h\in\{0, 1, \ldots, s-1\}$. Indeed, if $h\in\{s, 1, \ldots, 2s-1\}$, then we can write $h =  h'+s$, with $h' = 0, \ldots, s - 1$, so that the equation above can be written as
	$$
	\sum_{j=0\atop j \equiv h' + s  \ (2s)}^m a_j = \sum_{j=0\atop j \equiv h' \ (2s) }^m a_j.
	$$
	This returns the same equation because $j \equiv h' + s \equiv h'-s \ (2s)$. 
\end{proof}

 \section{Proof of the main results}\label{proofs}

\begin{proof}[Proof of Theorem \ref{T-conjecture}] 
	 Assume first that  $N$ is square-free, that is $\dsize
		N:=\prod_{r=1\atop p_r\in \Primes}^n p_r,\ \hbox{with}\ p_1<p_2 <\ldots < p_n.
		$   
	 
Let $s\in\N$ such that
	$ s\varphi(N)=s\prod_{r=1}^n (p_r-1)\le m =\deg {\cal P}(z)  $. We  show that
	${\cal P}(z) $ is divisible by  $\Phi_{N} (z^{s})$ if and only if 
	\begin{equation}\label{e-conj-p>2square-free}
		\sum_{d|N}(-1)^{\omega(d)}\sum_{j\equiv_{sN}\ h-s\sum_{p|d}\frac{N}{p}} a_j=0
	\end{equation}
	for every $h\in  \{0, 1,\ldots, sN-1\}$.
	
	To this end, we proceed by induction on $n=\omega(N)$, where  Lemma \ref{T-generaldivisibility} gives the  base case $n=1$. Thus, let us assume that   \eqref{e-conj-p>2square-free} is true for $n\ge 1$ and prove   that, given a prime $p_{n+1}>p_n$ such that 	$s\prod_{r=1}^{n+1} (p_r-1)\le m$, the polynomial ${\cal P}(z) $ is divisible by  $\Phi_{p_1\cdots p_n p_{n+1}} (z^{s})=\Phi_{Np_{n+1}} (z^{s})$ if and only if 
	\begin{equation}\label{e-conj-p>2n+1}
		\sum_{d|N p_{n+1}}(-1)^{\omega(d)}\sum_{j\equiv_{sN p_{n+1}}\ h-s\sum_{p|d}\frac{N p_{n+1}}{p}} a_j=0
	\end{equation}
	for every $h\in \{0, 1,\ldots, s N p_{n+1}-1\}$.
	For this purpose,  
	by using Lemma \ref{cycloformulas}  we 
	write 
	\begin{align*}
		\frac{{\cal P}(z)}{\Phi_{p_1\cdots p_n p_{n+1}}(z^{s}) } &= \frac{{\cal P}(z)}{\Phi_{Np_{n+1}} (z^{s})}= \frac{{\cal P}(z)}{\Phi_{N} (z^{sp_{n+1}})} {\Phi_{N} (z^{s  })}\\
		&=\frac{{\cal P}(z)}{\Phi_{N} (z^{sp_{n+1}})}	\prod_{d\vert N} (z^{sN/d}-1)^{\mu(d)}
		\\
		&=\frac{(z^{sN}-1){\cal P}(z)}{\Phi_{N } (z^{sp_{n+1}})}	
			 \prod_{d\vert N\atop{d>1}} (z^{sN/d}-1)^{\mu(d)},
	\end{align*}
	where recall that $\mu(1)=1$.
	
	Now, let us show that if $d>1$ divides $N$, then the polynomials
	$z^{sN/d} - 1$ and $\Phi_N(z^{s p_{n+1}})$  are coprime, i.e.,  they have no common roots. 
	
	To this end, first note that
	the roots of $z^{sN/d} - 1$ are 
	$$
	\zeta_{sN}^{jd}= e^{2\pi i jd / sN},\quad  j = 0, 1, \ldots, sN/d - 1.
	$$
Substituting $\zeta_{sN}^{jd}$ into $\Phi_N(z^{s p_{n+1}})$ gives
	$$
	\Phi_N(\zeta_{sN}^{jds p_{n+1}}) = \Phi_N\left(e^{2\pi i j d p_{n+1} / N}\right).
	$$
	Since $d =  (d, N) > 1$, we  see that  $e^{2\pi i j d p_{n+1} / N}$  cannot be a   primitive $N$th root of unity. Hence, $\zeta_{sN}^{jd}$ is not a root of $\Phi_N(z^{s p_{n+1}})$ for any $j = 0, 1, \ldots, sN/d - 1$.

	As a consequence,  $\Phi_{Np_{n+1} } (z^{s })$ divides ${\cal P}(z)$ if and only if  $\Phi_N(z^{s p_{n+1}})$ 
	divides  the polynomial
	$$ 
	(z^{sN}-1){\cal P}(z) = \sum_{j=0}^{m+sN} b_j z^{m+sN-j}, $$
	where $b_j= a_j-a_{j-sN}$. Recall that $a_i=0$ when $i>m$ and $i<0$. 
	
	By inductive assumption,  this is equivalent to
	
	\begin{equation}\label{1}
		S(h):= 	
		\sum_{d|N}(-1)^{\omega(d)}\sum_{j\equiv_{sNp_{n+1}}\ h-s\sum_{p|d}\frac{Np_{n+1}}{p}}   (a_j-a_{j-sN})=0
	\end{equation}
	for every  $h\in  \{0, 1,\ldots, sNp_{n+1}-1\}$.  
	
	For convenience, in the next formulas we let 
	$w=w( N, s, p_{n+1}):=sNp_{n+1}$ and $ f(d,h)= f(d, h, w):=h-s\sum_{p|d}\frac{Np_{n+1}}{p}=h-\sum_{p|d}\frac{w}{p}  $.
  
	Note that $j\equiv_w f(d,h)$ if and only if
	$ 
	j-sN \equiv_w
	f(dp_{n+1},h)$.

	Therefore, we see that
	\begin{align*}
		S(h)= &	\sum_{d|N}(-1)^{\omega(d)}\sum_{j\equiv_{w}\ f(d,h)} a_j\ + \ 
		\sum_{d|N}(-1)^{\omega(dp_{n+1})}\sum_{j\equiv_{w}\ f(dp_{p+1},h)} a_{j}\\
		= &	\sum_{d|N}(-1)^{\omega(d)}\sum_{j\equiv_{w}\ f(d,h)} a_j\ +  
		\sum_{d|Np_{n+1}\atop d\equiv 0\ (p_{n+1})}(-1)^{\omega(d)}\sum_{j\equiv_{w}\ f(d,h)} a_{j}\\
		=& 	\sum_{d|N p_{n+1}}(-1)^{\omega(d)}\sum_{j\equiv_{sN p_{n+1}}\ h-s\sum_{p|d}\frac{Np_{n+1}}{p}} a_j
	\end{align*}
	for every  $h\in  \{0, 1,\ldots, sNp_{n+1}-1\}$. 
	
	Hence,
	 \eqref{e-conj-p>2n+1} is an immediate consequence of \eqref{1}. 
	
	The equation \eqref{e-conj-p>2} for any $N$ follows straightforwardly from \eqref{e-conj-p>2square-free} after noticing that   \eqref{e-p-primes} yields
	$$
	\Phi_{N}(z)=\Phi_ {\gamma(N)}(z^{N/\gamma(n)}),
	$$
	where $
	\gamma(N)$ 
	is the  square-free kernel of $N$.

	When $N$ is even, an easy reformulation  of the  argument used in the proof  of Lemma \ref{T-generaldivisibility}  shows that it suffices to take $h\in \{0, 1,\ldots,N/2-1\}$.
	
	Theorem \ref{T-conjecture} is completely proved.
\end{proof}

\begin{Rem}\label{newrem}
	When $ \varphi(N)\le m < N  $, the  formula  \eqref{e-conj-p>2} can be simplified. Indeed, for any given  $h\in \{0, 1,\ldots,N-1\}$ and any divisor $ d $ of $ N $, one has that $  j \equiv h - \sum_{p \mid d} N/p\ \Mod N $ if and only if 
	$ 
	j = \{ h - \sum_{p \mid d} N/p \}_{N } + k N 
	$ 
	for some $ k \in \Z $. Recall that 
	$ 
	\{v\}_u =v- u [v/u ]
	$ for $ u, v \in \N $.
	
If, in addition, $ 0 \leq j \leq m < N  $,  then it must be $k=0$, i.e., 
	$ 
	j =  \{ h - \sum_{p \mid d} N/p \}_{N }
	$ 
	is the  only  solution of the congruence $  j \equiv h - \sum_{p \mid d} N/p\ \Mod N $. In this case, \eqref{e-conj-p>2} reduces to
	 
$$
		\sum_{d \mid \gamma(N)} (-1)^{\omega(d)} \, a_{\{ h - \sum_{p \mid d} N/p\}_{N }} = 0, \quad h\in \{0, 1,\ldots,N-1\}.
$$
	Clearly, these identities   are also satisfied by
	the coefficients of $\Phi_N(z)$.
\end{Rem}

 \begin{proof}[Proof of Corollary \ref{T-Ram-conj}] 
 	As already mentioned, the cyclotomic  polynomial 
 	$\Phi_n(z)$ divides the polynomial
 	\eqref{Toth1}, which can be written in the form \eqref{polynomial} as
 	$$
 	T_n(z)= \sum_{j=0}^{\tau} c_n(\tau-j) z^{\tau-j}-n 
 	 = \sum_{j=0}^{	\tau} a_j z^{\tau-j},
 	$$
 	where $\tau =\tau(n):= n-n/\gamma(n)$ and
 	\begin{equation}\label{a-coeff} a_j:=\begin{cases}  c_n(\tau-j) =c_n(j+n/\gamma(n))  &  \mbox{ if $0\le j <\tau$,}
 			\cr   c_n(0)-n =\varphi(n)-n & \mbox{ if $j=\tau$,}\cr
 			 0 & \mbox{ otherwise.}\cr
 	\end{cases}
 \end{equation}
 Therefore, since $\varphi(n)\le \tau= n-n/\gamma(n)<n$, from	Theorem \ref{T-conjecture}  (see also Remark \ref{newrem}) we get 
  \begin{equation}\label{e-conj-p>2rama}
 	\sum_{d|n}\mu(d)\sum_{j\equiv_{n}\ h-\sum_{p|d}\frac{n}{p}} a_j=
 \sum_{d|n}\mu(d)\, a_{\{h-\sum_{p|d}\frac{n}{p}\}_{n}}=0,
 \end{equation}
 for every 
 $$
 h\in H=H(n):=
 \begin{cases} \{0, \ldots, n/2-1\} & \mbox{if $n$ is even},\cr
 	\{0, \ldots, n-1\}   & \mbox{if $n$ is odd.}
 	\cr\end{cases}
 $$
Now,  given any divisor $d_1$ of $\gamma(n)$, let  
$h_1 \in H$  be such that $\{h_1-\sum_{p|d_1}\frac{n}{p}\}_{n}=\tau$, i.e.,
$h_1\equiv\ \sum_{p|d_1}n/p-n/\gamma(n)\ \Mod n$. Thus, 
\eqref{e-conj-p>2rama} becomes
$$
	\sum_{d|n\atop d\not=d_1}\mu(d) a_{\{\sum_{p|d_1}\frac{n}{p}-\sum_{p|d}\frac{n}{p}-\frac{n}{\gamma(n)}\}_n}+\mu(d_1)a_{\tau}=0,
$$
which gives \eqref{directidentity} because from \eqref{a-coeff} it follows that
$$
a_{\{\sum_{p|d_1}\frac{n}{p}-\sum_{p|d}\frac{n}{p}-\frac{n}{\gamma(n)}\}_n}=
c_n\Big(\sum_{p|d_1}\frac{n}{p}-\sum_{p|d}\frac{n}{p}\Big),
\qquad
a_{\tau}=\varphi(n)-n.
$$

In order to prove \eqref{e-sum-cn}, it suffices to note that 
if $h \in H$ is  such that $\tau=n-n/\gamma(n)>\{h- \sum_{p|d}n/p\}_n $, then in view of \eqref{a-coeff} we can write 
$a_{\{h-\sum_{p|d}\frac{n}{p}\}_{n}}=c_n(  h+\frac n{\gamma(n)}-\sum_{p\vert d} \frac np )$ in \eqref{e-conj-p>2rama}.
\smallskip

It remains to prove \eqref{e-sum-cnnonsquarefree}. For this purpose, assume that $\mu(n)=0$, i.e.,
$n/\gamma(n)\ge 2$, and observe that this yields $$
\Big(\sum_{p|d_1}n/p-n/\gamma(n),\ \sum_{p|d_1}n/p\Big)\cap\N\not=\emptyset
$$ 
for any divisor $d_1>1$ of $\gamma(n)$.
Since $h\in\big(\sum_{p|d_1}n/p-n/\gamma(n),\sum_{p|d_1}n/p\big)$ is equivalent to 
$\tau=n-n/\gamma(n)<h-\sum_{p|d_1}n/p+n<n$, it must be
$a_{\{h-\sum_{p|d_1}\frac{n}{p}\}_{n}}=0$. Hence, \eqref{e-sum-cnnonsquarefree} follows immediately from \eqref{e-conj-p>2rama}.

Corollary \ref{T-Ram-conj} is completely proved.\end{proof}

\begin{proof}[Proof of Theorem \ref{T-newcoeffformula}] Since   \eqref{e-primepower} of Lemma \ref{cycloformulas} yields  $	\Phi_{m}(z^p)=\Phi_{mp}(z)\Phi_m(z)$, we can write
	\begin{align*}
		\Phi_{m}(z^p)=&\sum _{j=0}^{\varphi(m)}a_m(j)z^{p\varphi(m)- pj}=
		\sum _{j'=0\atop j'\equiv 0\, (p)}^{p\varphi(m)}a_m(j'/p)z^{p\varphi(m)-j'} \\
		=&\Phi_{mp}(z)\Phi_m(z)=\left(\sum _{k=0}^{(p-1)\varphi(m)}a_{mp}(k)z^{(p-1)\varphi(m)-k}\right)\prod_{1\leq j\leq m\atop{(j,m)=1}} (z-\zeta_m^j),
	\end{align*}
	where recall that ${\deg }\Phi_{mp}=\varphi(pm)=(p-1)\varphi(m)$ for $p$ does not divide $m$.
	
	By applying  Lemma \ref{L-division-pol-3} to $\Phi_m(z^p)$, the formula \eqref{e-ck1} gives
	$$
	a_{mp}(k)=\sum_{j'=0\atop j'\equiv 0\, (p)}^ka_m(j'/p){\cal H}_{k-j'}(\vec\zeta_m^*)
	=
	\sum_{s=0}^{[k/p]} a_m(s){\cal H}_{k-sp}(\vec\zeta_m^*),
	$$
	for every $ k\in\{0,1,\ldots, (p-1)\varphi(m)\}$,
	as claimed in \eqref{newAM}.
\end{proof}  
\begin{proof}[Proof of Theorem \ref{T-coeff-pq}] By applying  Viete's  formula \eqref{e-Vieta} and Lemma \ref{L-rootsformulas} to
	$ 
	\Phi_{pq}(z)=\sum_{k=0}^{\varphi(pq)} a_{pq}(k)z^{\varphi(pq)-k},
	$  
	we get
	$$ 
	a_{pq}(k)= (-1)^k\E_k(\vec \zeta_{pq} ^*)=  {\cal H}_k(\vec \zeta_{pq} ^{**})
	$$
	for every $k=0,1,\ldots,\varphi(pq)=(p-1)(q-1)$. 
	
	Recall that
	$\vec \zeta_{pq} ^*$ denotes   the  vector of  all the primitive $pq$th roots of unity, which are given by $\zeta_p^k\zeta_q^s$, with $k=1,\ldots, p-1$ and $s=1,\ldots, q-1$. Thus,  $1,\zeta_p^k,\zeta_q^s$, with $k=1,\ldots, p-1$ and $s=1,\ldots, q-1$, are  the non-primitive $pq$th roots. Consequently,  
	we can write
	$$
	a_{pq}(k)={\cal H}_k(\vec \zeta_{pq} ^{**})={\cal H}_k
	(1, \zeta_p,\ldots, \zeta_p^{p-1}, \zeta_q,\ldots, \zeta_q^{q-1})
	$$
	for every $k=0,1,\ldots,\varphi(pq)=(p-1)(q-1)$.
	
	By using \eqref{e-id-2terms} of  Lemma \ref{L-id-homog} we  see that  
	\begin{equation}\label{coeffhom}
	a_{pq}(k)= {\cal H}_k (\vec\zeta_p ,\, \vec \zeta_q )-
	{\cal H}_{k-1}(\vec\zeta_p ,\, \vec \zeta_q ),
	\end{equation}
	where $\vec\zeta_p:=(1, \zeta_p,\ldots, \zeta_p^{p-1})$ and $\vec\zeta_q  :=(1, \zeta_q,\ldots, \zeta_q^{q-1})$.
	
	Now, let us apply \eqref{e-sigma-id 2} of Lemma \ref{L-id-homog} and the formula \eqref{e-sigma-id 3}  to write
	$$
	{\cal H}_k(\vec \zeta_p, \vec \zeta_q)= \sum_{s=0}^k  {\cal H}_{s} (\vec \zeta_p)  {\cal H}_{k-s}(\vec \zeta_q)= 
	\sum_{{s=0\atop {s\equiv 0\, \, (p)}}\atop s\equiv k\, \, (q)}^k  1
	$$
	Since $p\not=q$, by the Chinese Remainder Theorem the system 
	$$
	\begin{cases}  s\equiv 0\ \, (p)&  \mbox{}    
		\cr  s\equiv k\ \, (q)&   \mbox{}\cr\end{cases}
	$$
	is satisfied only by the integers $s\equiv kpv\ \, (pq)$, where $pv\equiv 1\ \, (q)$. However, being $k\le (p-1)(q-1)$, there is at most one integer $s\in[0,k]$ such that $s\equiv kpv\ \, (pq)$. Indeed, since
	it must be $0\le s=kpv+rpq\le k$ for some integer $r$, i.e.,
	$$
	-r\in\left[k\left(\frac{v}{q}-\frac{1}{pq}\right), k\frac{v}{q}\right]\cap\Z,
	$$
	it suffices to note that the length of this interval is
	$$
	\frac{k}{pq}\le\frac{(p-1)(q-1)}{pq}<1.
	$$
	Thus, ${\cal H}_k(\vec \zeta_p, \vec \zeta_q)=1$  
	for any $k=0,1,\ldots,(p-1)(q-1)$ such that
	$$
	\left[k\left(\frac{v}{q}-\frac{1}{pq}\right), k\frac{v}{q}\right]\cap\Z\not=\emptyset,
	$$
	and ${\cal H}_k(\vec \zeta_p, \vec \zeta_q)=0$  otherwise.
	Analogous conclusion holds for ${\cal H}_{k-1}(\vec \zeta_p, \vec \zeta_q)$. Hence, in view of \eqref{coeffhom} it must be $	a_{pq}(k)\in\{-1,0,\ 1\}$.
\end{proof}

\section{Appendix}\label{appendix}

\subsection{Alternate proof of (1.6)	}

For most of the considerations in this section   it is   tacitly assumed that we are dealing with square-free integers  in the support of  the $\mu$ function.  For example, we  freely use without explicit mention the fact that if $q$ is square-free, then $(d,q/d)=1$ for all $d|q$, so that $g(q)=g(d)g(q/d)$ for any multiplicative arithmetic function $g$ involved here. 

\medskip
First, notice that for any $d|\gamma(n)$ one has
$\big(n,\sum_{p|d}n/p\big)=n/d$. 

If $d_1=1$, then $\mu(d_1)=1$ and $\sum_{p|d_1}n/p=0$. Thus,
by H\"older's identity (see (4) of Lemma \ref{rama}) and \eqref{squarefree} we have that
$$			
\sum_{d|\gamma(n)\atop d\not=1}\mu(d)c_n\Big(-\sum_{p|d}\frac{n}{p}\Big)=
\varphi(n)
\sum_{d|n\atop d\not=1}\frac{\mu(d)^2}{\varphi(d)}=	\varphi(n)\Big(\frac{n}{\varphi(n)}-1\Big)		
=n-\varphi(n),
$$
which is \eqref{directidentity} for $d_1=1$.

Let us assume that $d_1>1$ and write
\begin{align}		
	\sum_{d|\gamma(n)\atop d\not=d_1}\mu(d)c_n\Big(	\sum_{p|d_1}\frac{n}{p}-\sum_{p|d}\frac{n}{p}\Big)&=c_n\Big(\sum_{p|d_1}\frac{n}{p}\Big)+
	\sum_{d|\gamma(n)\atop d\not=1,d_1}\mu(d)c_n\Big(\sum_{p|d_1}\frac{n}{p}-\sum_{p|d}\frac{n}{p}\Big)\nonumber\\
	&=\varphi(n)\frac{\mu(d_1)}{\varphi(d_1)}+
	\sum_{d|\gamma(n)\atop d\not=1,d_1}\mu(d)c_n\Big(\sum_{p|d_1}\frac{n}{p}-\sum_{p|d}\frac{n}{p}\Big).\label{step1}
\end{align}
For the sum on the right-hand side of \eqref{step1} we see that
\begin{align*}		
	\sum_{d|\gamma(n)\atop d\not=1,d_1}\mu(d)c_n\Big(\sum_{p|d_1}\frac{n}{p}-&\sum_{p|d}\frac{n}{p}\Big)\\
	&=\sum_{t|d_1}\sum_{{d|\gamma(n)\atop d\not=1,d_1}\atop (d,d_1)=t}\mu(d)c_n\Big(\sum_{p|\frac{d_1}{t}}\frac{n}{p}-\sum_{p|\frac{d}{t}}\frac{n}{p}\Big)\nonumber\\
	&=\sum_{t|d_1}\mu(t)\sum_{{d|\frac{\gamma(n)}{t}\atop dt\not=1,d_1}\atop (d,d_1/t)=1}\mu(d)c_n\Big(\sum_{p|\frac{d_1}{t}}\frac{n}{p}-\sum_{p|d}\frac{n}{p}\Big).
\end{align*}
Observe that the condition $dt=d_1$, with $(d,d_1/t)=1$, is satisfied if and only if $t=d_1$ and $d=1$.  Therefore, we can write
\begin{align}\label{step2}
	\sum_{d|\gamma(n)\atop d\not=1,d_1}\mu(d)c_n\Big(\sum_{p|d_1}\frac{n}{p}-\sum_{p|d}\frac{n}{p}\Big)
	&=\sum_{{d|\gamma(n)\atop (d,d_1)=1}\atop d>1}\mu(d)c_n\Big(\sum_{p|d_1}\frac{n}{p}-\sum_{p|d}\frac{n}{p}\Big)+\nonumber\\
	&\qquad \mu(d_1)\sum_{d|\frac{\gamma(n)}{d_1}\atop d>1}\mu(d)c_n\Big(\sum_{p|d}\frac{n}{p}\Big)+\nonumber\\
	&\quad \sum_{t|d_1\atop 1< t<d_1}\mu(t)\sum_{d|\frac{\gamma(n)}{t}\atop (d,d_1/t)=1}\mu(d)c_n\Big(\sum_{p|\frac{d_1}{t}}\frac{n}{p}-\sum_{p|d}\frac{n}{p}\Big).
\end{align}
As before, we have that
\begin{equation}\label{step3}		
	\sum_{d|\frac{\gamma(n)}{d_1}\atop d>1}\mu(d)c_n\Big(\sum_{p|d}\frac{n}{p}\Big)=
	\varphi(n)
	\sum_{d|\frac{\gamma(n)}{d_1}\atop d>1}\frac{\mu(d)^2}{\varphi(d)}
	=	\frac{\gamma(n)}{d_1}
	\frac{\varphi(n)}{\varphi(\frac{\gamma(n)}{d_1})}-\varphi(n).
\end{equation}
Further,  since 
$$
\Big(n,\sum_{p|d_1}\frac{n}{p}-\sum_{p|d}\frac{n}{p}\Big)=\frac{n}{dd_1},
$$ 
by H\"older's identity and the condition $(d,d_1)=1$ one has
$$
c_n\Big(\sum_{p|d_1}\frac{n}{p}-\sum_{p|d}\frac{n}{p}\Big)=\varphi(n)\frac{\mu(d)\mu(d_1)}{\varphi(d)\varphi(d_1)}.
$$
Thus, the first sum on the right-hand side of \eqref{step2} becomes
\begin{align*}
	\sum_{{d|\gamma(n)\atop (d,d_1)=1}\atop d>1}\mu(d)c_n\Big(\sum_{p|d_1}\frac{n}{p}-\sum_{p|d}\frac{n}{p}\Big)
	&=\varphi(n)\frac{\mu(d_1)}{\varphi(d_1)}
	\sum_{{d|\gamma(n)\atop (d,d_1)=1}\atop d>1}\frac{\mu(d)^2}{\varphi(d)}
	\\
	&=\varphi(n)\frac{\mu(d_1)}{\varphi(d_1)}\left(
	\sum_{d|\gamma(n)\atop (d,d_1)=1}\frac{\mu(d)^2}{\varphi(d)}-1\right).
\end{align*}
By using Lemma \ref{eulerandmebius} we see that
\begin{align*} 
	\sum_{d|\gamma(n)\atop (d,d_1)=1}\frac{\mu(d)^2}{\varphi(d)}&=
	\sum_{d|\gamma(n)}\frac{\mu(d)^2}{\varphi(d)}
	\sum_{r| (d,d_1)}\mu(r)=\sum_{r| d_1}\mu(r)
	\sum_{d|\gamma(n)\atop d\equiv 0\ (r)}\frac{\mu(d)^2}{\varphi(d)}\\
	&=\sum_{r| d_1}\frac{\mu(r)}{\varphi(r)}
	\sum_{d|\frac{\gamma(n)}{r}}\frac{\mu(d)^2}{\varphi(d)}
	=\sum_{r| d_1}\frac{\mu(r)}{\varphi(r)}
	\frac{\gamma(n)/r}{\varphi(\gamma(n)/r)}\\
	&=	\frac{\gamma(n)}{\varphi\big(\gamma(n)\big)}\frac{\varphi(d_1)}{d_1}.
\end{align*}
Therefore,
\begin{align}\label{step4}
	\sum_{{d|\gamma(n)\atop (d,d_1)=1}\atop d>1}\mu(d)c_n\Big(\sum_{p|d_1}\frac{n}{p}-\sum_{p|d}\frac{n}{p}\Big)&=
	\varphi(n)\frac{\mu(d_1)}{d_1}\frac{\gamma(n)}{\varphi\big(\gamma(n)\big)}-
	\varphi(n)\frac{\mu(d_1)}{\varphi(d_1)}\nonumber\\
	&=
	n\frac{\mu(d_1)}{d_1}-
	\varphi(n)\frac{\mu(d_1)}{\varphi(d_1)},
\end{align}
after recalling that $\frac{\gamma(n)}{\varphi\big(\gamma(n)\big)}=\frac{n}{\varphi(n)}$ by \eqref{totient} of Lemma \ref{eulerandmebius}.

Analogously, for the third sum on the right-had side of \eqref{step2} one has
\begin{align*} 
	\sum_{t|d_1\atop 1< t<d_1}\mu(t)\sum_{d|\frac{\gamma(n)}{t}\atop (d,d_1/t)=1}&\mu(d)c_n\Big(\sum_{p|\frac{d_1}{t}}\frac{n}{p}-\sum_{p|d}\frac{n}{p}\Big)\\
	&=\mu(d_1)\varphi(n)
	\sum_{t|d_1\atop 1< t<d_1}\frac{1}{\varphi(d_1/t)}\sum_{d|\frac{\gamma(n)}{t}\atop (d,d_1/t)=1}\frac{\mu(d)^2}{\varphi(d)}\\
	&=\mu(d_1)\varphi(n)
	\sum_{t|d_1\atop 1< t<d_1}\frac{1}{\varphi(d_1/t)}\sum_{d|\frac{\gamma(n)}{t}}\frac{\mu(d)^2}{\varphi(d)}\sum_{r| (d,d_1/t)}\mu(r)\\
	&= \mu(d_1)\varphi(n)\sum_{t|d_1\atop 1< t<d_1}\frac{1}{\varphi(d_1/t)}
	\sum_{r|\frac{d_1}{t}}\frac{\mu(r)}{\varphi(r)}
	\sum_{d|\frac{\gamma(n)}{rt}}\frac{\mu(d)^2}{\varphi(d)}.
\end{align*}
Now, observe that 
\begin{align*}
	\sum_{t|d_1\atop 1< t<d_1}\frac{1}{\varphi(d_1/t)}&
	\sum_{r|\frac{d_1}{t}}\frac{\mu(r)}{\varphi(r)}
	\sum_{d|\frac{\gamma(n)}{rt}}\frac{\mu(d)^2}{\varphi(d)}\\
	&=\gamma(n)
	\sum_{t|d_1\atop 1< t<d_1}\frac{1}{t\varphi(d_1/t)}
	\sum_{r|\frac{d_1}{t}}\frac{\mu(r)}{r\varphi(r)\varphi\big(\gamma(n)/rt\big)}\\
	&=\gamma(n)
	\sum_{t|d_1\atop 1< t<d_1}\frac{1}{t\varphi(d_1/t)\varphi\big(\gamma(n)/t\big)}
	\sum_{r|\frac{d_1}{t}}
	\frac{\mu(r)}{r}\\
	&=\frac{\gamma(n)}{d_1}
	\sum_{t|d_1\atop 1< t<d_1}\frac{1}{\varphi\big(\gamma(n)/t\big)}
	\\
	&=\frac{\gamma(n)}{d_1}\left(
	\sum_{t|d_1}\frac{1}{\varphi\big(\gamma(n)/t\big)}-
	\frac{1}{\varphi\big(\gamma(n)\big)}
	-\frac{1}{\varphi\big(\gamma(n)/d_1\big)}
	\right)
	\\
	&=\frac{\gamma(n)}{d_1\varphi\big(\gamma(n)/d_1\big)}
	\sum_{t|d_1}\frac{1}{\varphi(d_1/t)}-
	\frac{\gamma(n)}{d_1\varphi\big(\gamma(n)\big)}
	-\frac{\gamma(n)}{d_1\varphi\big(\gamma(n)/d_1\big)}
	\\
	&=\frac{\gamma(n)}{\varphi\big(\gamma(n)\big)}
	-\frac{\gamma(n)}{d_1\varphi\big(\gamma(n)\big)}
	-\frac{\gamma(n)}{d_1\varphi\big(\gamma(n)/d_1\big)}.
\end{align*}
Here we have used the fact that for $\mu(d_1)\not=0$ one has
$$
\sum_{t|d_1}\frac{1}{\varphi(d_1/t)}= \sum_{t|d_1}\frac{\mu(t)^2}{\varphi(t)}.
$$
Thus,
\begin{align*}
	\sum_{t|d_1\atop 1< t<d_1}\mu(t)\sum_{d|\frac{\gamma(n)}{t}\atop (d,d_1/t)=1}&\mu(d)c_n\Big(\sum_{p|\frac{d_1}{t}}\frac{n}{p}-\sum_{p|d}\frac{n}{p}\Big)\\
	&=
	\mu(d_1)n-\frac{\mu(d_1)}{d_1}n-
	\varphi(n)\frac{\mu(d_1)}{d_1}\frac{\gamma(n)}{\varphi\big(\gamma(n)/d_1\big)}.
\end{align*}	
Together with \eqref{step1}-\eqref{step4}, this yields \eqref{directidentity} when $d_1\not=1$.	

\subsection{
Corollary \ref{T-Ram-conj} when \boldmath{$\omega(n)=2$}.}

Here we explicitly exhibit the case $\omega(n)=2$ of the Corollary \ref{T-Ram-conj}.  

Let  $n:=p_1^{v_1}p_2^{v_2}$, with  $v_1,v_2\in\N$, $p_1,p_2\in\Primes$ such that $p_1<p_2$.

The equation \eqref{directidentity}		 of
Corollary \ref{T-Ram-conj}  for $n=p_1^{v_1}p_2^{v_2}$  gives the identity
	$$
		c_n\Big(\frac{n}{p_2}\pm \frac{n}{p_1}\Big)=
		c_n\Big(\frac{n}{p_1}\Big)+c_n\Big(\frac{n}{p_2}\Big)+n-\varphi(n).
	$$
By using (8) of Lemma \ref{rama} for $p_1=2$, this formula reduces to
	$$
	2c_n(n/p_2)+n=2\varphi(n).
	$$
	The equation \eqref{e-sum-cn} of
Corollary \ref{T-Ram-conj}  for $n=p_1^{v_1}p_2^{v_2}$ says that
if 	
	$$
h\in	H:=
	\begin{cases} 		\{0, \ldots,n/2-1\}   & \mbox{if $p_1=2$,}\cr\{0, \ldots, n-1\} & \mbox{if $p_1>2$},
		\cr\end{cases}
	$$ 
	is  such that $\{h- \sum_{p|d}n/p\}_n < n-n/\gamma(n)$ 
for every $d\in\{1,p_1,p_2,p_1p_2\}$,   then
$$
\sum_{d\in\{1,p_1,p_2,p_1p_2\}}(-1)^{\omega(d)}c_{n}\Big(h+\frac{n}{p_1p_2} -\sum_{p|d}\frac{n}{p}\ \Big)=0.
$$
This is equivalent to
	\begin{align*}
		c_{n}\Big(h+\frac{n}{p_1p_2} \Big)+&
		c_{n}\Big(h+\frac{n}{p_1p_2}-\frac{n}{p_1}-\frac{n}{p_2}\ \Big)=\nonumber\\ 
		&
		c_{n}\Big(h+\frac{n}{p_1p_2} -\frac{n}{p_1}\ \Big)+c_{n}\Big(h+\frac{n}{p_1p_2}-\frac{n}{p_2}\ \Big)
	\end{align*}
	for every integer
	$$
	h\in\Big[0,\frac{n}{p_2}-\frac{n}{p_1p_2}\Big)\cup
	\Big[\frac{n}{p_2},\frac{n}{p_1}-\frac{n}{p_1p_2}\Big)\cup
	\Big[\frac{n}{p_1},\frac{n}{p_1}+\frac{n}{p_2}-\frac{n}{p_1p_2}\Big)\cup
	\Big[\frac{n}{p_1}+\frac{n}{p_2}, n-\frac{n}{p_1p_2}\Big).
	$$
	For $p_1=2$ this identity reduces to
	$$
	c_n\Big(h+\frac{n}{2p_2}\Big)=
	c_n\Big(h-\frac{n}{2p_2}\Big)
	$$
	for every integer $h\in\Big[0,\frac{n}{2p_2}\Big)\cup\Big[\frac{n}{p_2},\frac{n}{2}-\frac{n}{2p_2}\Big)$.
	
Finally, assuming that $\mu(n)=0$, i.e., $v_1v_2\ge 2$, the equation \eqref{e-sum-cnnonsquarefree} gives the following identities.	

\begin{itemize}
\item[1.] $c_n\Big(h+\frac{n}{p_1p_2}\Big)+c_n\Big(h-\frac{n}{p_1}-\frac{n}{p_2}+\frac{n}{p_1p_2}\Big)=
	c_n\Big(h-\frac{n}{p_1}+\frac{n}{p_1p_2}\Big)$	
	
	for every  $h\in\Big(\frac{n}{p_2}-\frac{n}{p_1p_2}, \frac{n}{p_2}\Big)\cap\N$.
	
	\item[2.] 	$c_n\Big(h+\frac{n}{p_1p_2}\Big)+c_n\Big(h-\frac{n}{p_1}-\frac{n}{p_2}+\frac{n}{p_1p_2}\Big)=
	c_n\Big(h-\frac{n}{p_2}+\frac{n}{p_1p_2}\Big)$
	
	for every  $h\in \Big(\frac{n}{p_1}-\frac{n}{p_1p_2}, \frac{n}{p_1}\Big)\cap\N$.
	
	\item[3.] $c_n\Big(h+\frac{n}{p_1p_2}\Big)=
	c_n\Big(h-\frac{n}{p_1}+\frac{n}{p_1p_2}\Big)+c_n\Big(h-\frac{n}{p_2}+\frac{n}{p_1p_2}\Big)$

for every $h\in  \Big(\frac{n}{p_1}+\frac{n}{p_2}-\frac{n}{p_1p_2}, \frac{n}{p_1}+\frac{n}{p_2}\Big)\cap\N$.

\item[4.] $c_n\Big(h-\frac{n}{p_1}-\frac{n}{p_2}+\frac{n}{p_1p_2}\Big)=
c_n\Big(h-\frac{n}{p_1}+\frac{n}{p_1p_2}\Big)+c_n\Big(h-\frac{n}{p_2}+
\frac{n}{p_1p_2}\Big)$

for every $h\in \Big(n-\frac{n}{p_1p_2},\ n-1\Big]\cap\N$.
\end{itemize}

For $p_1=2$ these formulas reduce to 
$$
c_n\Big(h+\frac{n}{2p_2}\Big)=2
c_n\Big(h-\frac{n}{2p_2}\Big)
$$
for every  $h\in\Big(\frac{n}{2}-\frac{n}{2p_2},\frac{n}{2}\Big)\cap\N$.

\end{document}